\documentstyle[12pt]{article}   

 \title{
Actions of the Dipper-Donkin quantization $GL_2$ on the Clifford 
algebra ${\it C}(1,3)$. }    
	
\author{
Suemi Rodr\'{\i}guez-Romo $^{\ast}$\\
Centro de Investigaciones Te\'oricas,\\ 
Universidad Nacional Aut\'onoma de M\'exico, 
Campus Cuautitl\'an, \\
Apdo. Postal 142, Cuautitl\'an Izcalli, 
Edo. de M\'exico, \\
54740 M\'exico}         
\date{}  
\begin{document}   
\maketitle 
\renewcommand{\thefootnote}{\fnsymbol{footnote}}
\setcounter{footnote}{-1}
\footnote{$\hspace*{-6mm}^{\ast}$
e-mail: suemi$@$servidor.unam.mx }
\renewcommand{\thefootnote}{\arabic{footnote}}\

\baselineskip0.6cm

{\small {\bf Abstract. } Following the method already developed for
studying the actions of $GL_q(2,C)$ on the Clifford algebra ${\it C}(1,3)$ 
and
its quantum invariants \cite{aqg1}, we study the action on 
${\it C}(1,3)$ of the quantum $GL_2$ constructed by Dipper and 
Donkin \cite{Di}. We are able of proving that there exits only two 
non-equivalent cases of actions with nontrivial ``perturbation" \cite{aqg1}. 
The spaces of invariants are trivial in both cases.\

We also prove that each irreducible finite dimensional algebra representation
of the quantum $GL_2$, $q^m\neq 1$, is one dimensional.

By studying the cases with zero ``perturbation" we find that the cases
with nonzero ``perturbation" are the only ones with maximal possible 
dimension for the operator algebra $\Re$.
\newpage
\baselineskip0.6cm
\newtheorem{theorem}{Theorem}
\newtheorem{corollary}{Corollary}
\newtheorem{lemma}{Lemma}  
\section{Introduction.}           
In this paper we consider inner actions of the Dipper-Donkin 
quantization of $GL_2$ (see \cite{Di}) on the space-time Clifford 
algebra ${\it C}(1,3)$. The analogous problems for Manin quantization 
was considered in details in \cite{aqg1}. We prove that every irreducible
finite dimensional algebra representation of $GL_2$, $q\neq 1$, is one 
dimensional and therefore triangular. Using this fact 
we show that only two particular cases have nonzero 
{\it perturbation}. Actions with trivial {\it perturbation} are also studied. 
From this, some consequences are derived.\

This paper is organized as follows. In Section 2, we introduce
elementary notions. In Section 3 we prove
Theorem 1 which is fundamental to address our problem. This Theorem 
deals with $q$-spinor representations ($q^3\neq 1$, 
$q^4\neq 1$), corresponding to representations of $GL_2$ with non 
zero {\it perturbation}. Since each irreducible finite dimensional algebra 
representation of Dipper-Donkin $GL_2$ is one dimensional, we can use the 
method of \cite{aqg1} in complete generality for the classification of 
inner actions. Finally, in Section 4 we present the representations of 
$GL_2$ by Dipper-Donkin with nonzero {\it perturbation} and some 
remarkable features related with the zero {\it perturbation} cases.\ 
\section{Preliminary notions.}
The algebraic structure of Dipper-Donkin quantization $GL_2$ 
\cite{Di} is generated by four elements $c_{ij}$, 
$1\leq i,j\leq 2$ with relations which can be presented by the 
following diagram.
%\vspace*{1.5cm}
\newline
\hspace*{2cm}
\begin{picture}(0,80)
\put(115.29,2.39){\makebox(0,0)[cc]{$c_{21}$}}
\put(141,70.5){\makebox(0,0)[cc]{$d$}}
\put(168.32,2.39){\makebox(0,0)[cc]{$c_{22}$}}
\put(168.32,50.46){\makebox(0,0)[cc]{$c_{12}$}}
\put(115.29,50.46){\makebox(0,0)[cc]{$c_{11}$}}
\put(117.50,7){\vector(1,1){40}}
\put(115.5,7){\vector(0,1){40}}
\put(168.5,7){\vector(0,1){40}}
\put(115.7,7){\vector(1,2){28}}
\put(147,68){\vector(3,-2){18}}
\put(120.5,3){\line(1,0){40}}
\put(120.5,51.0){\line(1,0){40}}
\put(115.5,55){\line(5,3){18}}
\put(146,63){\line(2,-5){23}}
\put(123.5,42.0){\makebox(0,0)[cc]{$\cdot $}}
\put(131.5,34.0){\makebox(0,0)[cc]{$\cdot $}}
\put(139.50,26.0){\makebox(0,0)[cc]{$\cdot $}}
\put(147.5,18.0){\makebox(0,0)[cc]{$\cdot $}}
\put(156.5,10.0){\makebox(0,0)[cc]{$\cdot $}}
\end{picture}
\vspace*{.5cm}
\label{fig1}
\newline
\hspace*{6cm} Figure 1. $GL_2$\\

Here we denote by arrows the ``quantum spinors'' ( or 
generators of the quantum plane \cite{Man}) $xy=qyx$ by 
straigth line the ``classical spinors'' $xy=yx$ \cite{aqg1} 
and by dots a classical spinor with a nontrivial 
{\it perturbation} \cite{aqg1}, $xy-yx=p$ being
$p$=$(q-1)c_{12}c_{21}$.\

Here the quantum determinant $d=c_{11}c_{22}-c_{12}c_{21}$
 is noncentral and group-like. This, in contrast with Manin's 
approach \cite{Man}. In any Hopf algebra every group-like 
element is invertible, therefore the quantum $GL_2$ includes the
formal inverse $d^{-1}$.\

The coalgebra structure is defined in the standard way for all
quantizations and the antipode $S$ is given in reference [2].\

The Clifford algebra ${\it C}(1,3)$ is generated by the vector 
$\gamma_{\mu}$, $\mu=0,1,2,3$ with relations defined by the form 
$g_{\mu\nu}$=$diag(1, -1, -1, -1)$, as follows:
$$
\gamma_{\mu}\gamma_{\nu}=g_{\mu\nu}+\gamma_{\mu\nu},\;\;
\gamma_{\mu\nu}=-\gamma_{\nu\mu},
$$
$$
\gamma_{\rho}\gamma_{\mu\nu}=
g_{\rho\mu}\gamma_{\nu}-g_{\rho\nu}\gamma_{\mu}
+\gamma_{\rho\mu\nu},
$$
$$
\gamma_{\lambda}\gamma_{\mu\nu\rho}=g_{\lambda\mu}\gamma_{\nu\rho}-
g_{\lambda\nu}\gamma_{\mu\rho}+g_{\lambda\rho}\gamma_{\mu\nu}+
\gamma_{\lambda\mu\nu\rho}.
$$
This algebra is isomorphic to the algebra of the $4\times 4$ complex 
matrix and it has the basis of matrix units reported in reference
\cite {aqg1}, among others.\

An action of $GL_2$ on ${\it C}(1,3)$ is uniquely defined by actions of
$c_{ij}$ on the generators of ${\it C}(1,3)$\cite{Co}\cite{Sh};
\begin{equation}
c_{ij}\cdot \gamma_k=
f_{ijk}(\gamma_0, \gamma_1, \gamma_2, \gamma_3),
\label{ac}
\end{equation}
where $f_{ijk}$ are some noncommutative polynomials in four 
variables. If (\ref{ac}) defines an action of quantum group 
$GL_2$ on ${\it C}(1,3)$ and $\gamma'_0$, $\gamma'_1$, 
$\gamma'_2$, $\gamma'_3$ is another  system of generators of 
${\it C}(1,3)$, with the same relations, then the formula
\begin{equation}
c_{ij}*\gamma'_k=
f_{ijk}(\gamma'_0, \gamma'_1, \gamma'_2, \gamma'_3),
\label{ac*}
\end{equation}
with the same polynomials $f_{ijk}$, will also define an 
action of the quantum $GL_2$ on ${\it C}(1,3)$. Two actions of 
$GL_2$ on ${\it C}(1,3)$ are said to be 
equivalent if they can be presented as in (\ref{ac}) and (\ref{ac*}) with 
the same polynomials $f_{ijk}$}. It is easy to show that two actions
$\cdot$, $*$ are equivalent if and only if $c_{ij}*(uwu^{-1})$=
$u(c_{ij}\cdot w)u^{-1}$ for some invertible $u\in {\it C}(1,3)$
(see \cite{aqg1}, formula (7)).\

For every action $\cdot$ there exist an invertible matrix
$M=\left(
\begin{array}{cc}
m_{11} & m_{12} \\
m_{21} & m_{22}
\end{array}
\right)$
$\in {\it C}(1,3)_{2\times 2}$, such that
$$
c_{ij}\cdot v=\sum m_{ik}vm^*_{kj},
$$
where 
$
\left(
\begin{array}{cc}
m^*_{11} & m^*_{12} \\
m^*_{11} & m^*_{12}
\end{array}
\right)
$=$M^{-1}$ (see Skolem-Noether theorem for Hopf algebras
\cite{Ko}\cite{Mo}). The action $\cdot$ is called inner if the map
$c_{ij}\rightarrow m_{ij}$ defines an algebra homorphism
$\varphi: GL_2\rightarrow {\it C}(1,3)$. Since the algebra
${\it C}(1,3)$ is isomorphic to the algebra of $4\times 4$
matrices, the homorphism ${\it C}(1,3)$ defines (and is
defined by) a four dimensional module over (the algebraic
structure of) $GL_2$, or, equivalently, a four dimensional
representation of $GL_2$.\

If $\varphi(c_{12}c_{21})$=$0$,
then by definition in Figure 1 the representation $\varphi$ is
defined for an essentially more simple structure, generated by
two commuting ``quantum spinors" $(c_{21},c_{11})$ and 
$(c_{22},c_{12})$. Firstable we focus our attention on the
case when $\varphi(c_{12}c_{21})$$\neq 0$ and in this case we
say that the inner action defined by $\varphi$ has nonzero
{\it perturbation}.\

If we add a formal inverse $c^{-1}_{11}$, then the algebraic structure 
of Dipper-Donkin quantization
$GL_2$ is generated by the elements in the following diagram.\

\hspace*{2cm}
\begin{picture}(0,80)
\put(115.29,2.39){\makebox(0,0)[cc]{$c_{21}$}}
\put(168.32,2.39){$a_{22}=c_{11}^{-1}d$}
\put(168.32,50.46){$a_{12}=c_{11}^{-1}c_{12}$}
\put(115.29,50.46){\makebox(0,0)[cc]{$c_{11}$}}
\put(117.50,7){\line(1,1){40}}
\put(122,46){\line(1,-1){40}}
\put(115.5,7){\vector(0,1){40}}
\put(168.5,7){\vector(0,1){40}}
\put(120.5,3){\line(1,0){40}}
\put(120.5,51.0){\line(1,0){40}}
\end{picture}
\vspace*{.5cm}
\label{fig1}
\newline
\hspace*{6cm} Figure 2. $GL_2$\\

From here, it follows straightforward that, up to invertibility 
of $c_{11}$, the algebraic structure of
$GL_2$ can be considered like a tensor product 
${\aleph}\otimes {\aleph}$ where ${\aleph}$ is
the quantum plane.\

In the next Section we study $q$-spinors suitable of being used 
to represent the quantum $GL_2$. Concretely speaking we consider 
in details the following triangle.\\
\hspace*{2cm}
\begin{picture}(0,80)
\put(115.29,2.39){\makebox(0,0)[cc]{$C$}}
\put(168.32,50.46){$B$}
\put(115.29,50.46){\makebox(0,0)[cc]{$A$}}
\put(117.50,7){\vector(1,1){40}}
\put(115.5,7){\vector(0,1){40}}
\put(120.5,51.0){\vector(1,0){40}}
\end{picture}
\vspace*{.5cm}
\label{fig1}
\newline
\hspace*{6cm} Figure 3.\\ 

corresponding to\

\hspace*{2cm}
\begin{picture}(0,80)
\put(115.29,2.39){\makebox(0,0)[cc]{$c_{21}$}}
\put(168.32,50.46){$c_{12}$}
\put(115.29,50.46){\makebox(0,0)[cc]{$d$}}
\put(117.50,7){\vector(1,1){40}}
\put(115.5,7){\vector(0,1){40}}
\put(120.5,51.0){\vector(1,0){40}}
\end{picture}
\vspace*{.5cm}
\label{fig1}
\newline
\hspace*{6cm} Figure 4.\\ 

in Figure 1.\

We say that the representation of the $q$-spinor $xy$=$qyx$,
$x\rightarrow A$, $y\rightarrow B$ is {\it admissible} if 
there exists $C$ such that $x\rightarrow C$, 
$y\rightarrow B$ and $x\rightarrow C$, $y\rightarrow A$ are also a 
representation of $q$-spinor with $CB\neq 0$. In other words it
means that $d\rightarrow A$, $c_{12}\rightarrow B$, 
$c_{21}\rightarrow C$ is a representation of the subalgebra of $GL_2$, 
generated by $d$, $c_{12}$, $c_{21}$ with $CB\neq 0$.
\section{$q$-spinor representations.}           
Let $(x,y)$ be a $q$-spinor, $xy=qyx$.  
If $x\rightarrow A$, $y\rightarrow B$ is its 
representation by $4\times 4$ matrices over complex 
numbers, then for every invertible $4\times 4$ matrix 
$u$ and nonzero number $\alpha$, the map 
$x\rightarrow uAu^{-1}\alpha$, 
$y\rightarrow uBu^{-1}\alpha$ also defines a 
representation of the $q$-spinor. Following, \cite{aqg1}, 
we consider this two representations as {\it equivalent} 
ones. Thus, under investigation of representations of a 
$q$-spinor, we can suppose that the matrix $A$ has a 
Jordan Normal form and one of it's eigenvalues is equal 
to 1 (if $A\neq 0)$.\  

For a given matrix $A$, we denote by $B(A)$ the linear space of
all matrices $B$, such that $AB$=$qBA$ and by $B'(A)$ the set of 
all matrices $B'$ such that $B'A$=$qAB'$.
\begin{theorem}$\!\!\!.$
Every admissible representation of the $q$-spinor 
$(q^3, q^4\neq 1)$ \cite{Mi}\cite{Sm}
by $4\times 4$ complex matrices, $x\rightarrow A$, 
$y\rightarrow B$, such that $A$ is an invertible matrix is equivalent 
to one of the following representations.  
\begin{eqnarray}
& 1. & A=diag(q^2, q, q, 1), \hspace{3cm}B=qe_{13}-\mu e_{24}\label{case1}\\
&    & \hspace*{6.7cm}B'=e_{43}-\mu e_{21}
\nonumber\\ 
& 2. & A=diag(q^2, q, q, 1),\hspace{3cm} B=qe_{12}+\mu e_{34}\label{case2}\\
&    & \hspace*{6.7cm}B'=e_{42}+\mu e_{31}
\nonumber\\
& 3. & A=diag\left(\left(
\begin{array}{cc}
q & 1\\
0 & q\\
\end{array}
\right), q^2, 1\right), \hspace*{1.3cm}B_1=e_{14}\;;\;B_2=e_{32}\label{case3}\\
&    & \hspace*{6.6cm}B'_1=e_{13}\;;\;B'_2=e_{42} \nonumber\\ \nonumber
\end{eqnarray}
\end{theorem}
{\it Proof}. If $B(A)^2\neq 0$; then by Theorem 1 \cite{aqg1}, 
we have seven different possibilities for $A$. Direct 
calculations show that only in the second case there exist 
representations with nonzero perturbation, these are (\ref{case1}) and
(\ref{case2}) described in the theorem.\

Let us now study the case $B(A)^2$=0. We assume that the matrix $A$ 
has a Jordan Normal form and one of 
its eigenvalues is equal to 1. By lemma 2 \cite{aqg1} the matrix $A$ 
cannot be a simplest Jordan Normal matrix; i.e. it has more than one 
block.\

If $A$=$diag(\alpha_1, \alpha_2, \alpha_3, 1)$ is a diagonal matrix 
then $B'(A)$ evidentely coincides with the space of transposed matrices 
$B(A)^T$. By Lemma 4 \cite{aqg1} 
the space $B(A)$ as well as $B(A)'$ are generated by matrix units 
and, by formula (20) \cite{aqg1}, $e_{ij}\in B(A)$ if and only if 
$\alpha_i$=$q\alpha_j$, (21) in \cite {aqg1}. Thus we have two main 
cases with $B(A)^2$=$0$; $\alpha_1$=$q\alpha_2$, $\alpha_3$=$q$; 
$e_{12}$, 
$e_{34}$$\in B(A)$ and $\alpha_1$=$\alpha_2$=$\alpha_3$=$q$;
$e_{12}$, $e_{13}$, $e_{14}$$\in B(A)$ while the others can be 
obtained from these by changing the numerations of indeces. In both
cases $\left(B(A)\cdot B(A)^T\right)\cap $
$\left(B(A)^T\cdot B(A)\right)$=$0$ and so there is no  admissible
representations.\

Let $A$ be of the form\newline
\begin{eqnarray}
A=
\left(
\begin{array}{cc}
a & 0 \\
0 & b
\end{array}
\right)
\label{Ama}
\end{eqnarray}
where $a$, $b$ are either invertible $2\times 2$ matrices in
Jordan Normal form or $a$ is an invertible simplest Normal Jordan 
$3\times 3$ matrix and $b$ is a nonzero complex number (and therefore
we can suppose that $b$=$1$).\

If $B'$=
$
\left(
\begin{array}{cc}
\alpha' & \beta' \\ 
\gamma' & \delta'
\end{array}
\right)
$
is a nonzero matrix from $B'(A)$ then by formula (23) \cite {aqg1} 
changing $q$ by $q^{-1}$ we have that\newline
\begin{eqnarray}
a\alpha'=q^{-1}\alpha' a &\;\;\;& a\beta'=q^{-1}\beta' b \\ \label{s}
b\gamma' =q^{-1}\gamma'a &\;\;\;& b\delta'=q^{-1}\delta'b. \label{su}
\end{eqnarray}
At first, let us consider when $a$ is a $3\times 3$ 
matrix. In \cite{aqg1} we can see that in this case 
there exists only two possibilities with $B(A)$$\neq 0$;\newline
$$
A=
\left(
\begin{array}{cccc}
q^{-1} & 1 & 0 & 0 \\
0 & q^{-1} & 1 & 0 \\
0 & 0 & q^{-1} & 0 \\
0 & 0 & 0 & 1
\end{array}
\right),\;\;
B(A)=Ce_{43}
$$
and
$$
A=
\left(
\begin{array}{cccc}
q & 1 & 0 & 0 \\
0 & q & 1 & 0 \\
0 & 0 & q & 0 \\
0 & 0 & 0 & 1
\end{array}
\right),\;\;
B(A)=Ce_{14}.
$$
In the first case, we have $B'(A)$=$Ce_{14}$ and in the second 
$B'(A)$=$Ce_{43}$. Thus the equality $Ce_{43}\cdot Ce_{14}$=$0$ 
shows that in both cases either $B(A)B'(A)$=$0$ or $B'(A)B(A)$
$=0$ and there exits no admissible representation.\

Consider now the case when $a, b, \alpha, \beta, \gamma, 
\delta$ are $2\times 2$ matrices. Here, we have defined 
$B=
\left(
\begin{array}{cc}
\alpha & \beta\\
\gamma & \delta\\
\end{array}
\right)$.\

Let us start with the case when both matrices $a$ and $b$ 
have a simplest Jordan Normal Form i.e.
\begin{equation}
a=
\left(
\begin{array}{cc}
\epsilon & 1\\
0	 & \epsilon
\end{array}
\right),\;\;\;
b=
\left(
\begin{array}{cc}
1& 1\\
0& 1
\end{array}
\right)
\label{sue}
\end{equation}
(recall that we suppose that one of the eigenvalues of 
$A$=diag$(a,b)$ is
equal to 1).\

We know that $[A,B]_q$=$0$  (from this follows that $\alpha$=
$\delta=0)$ and $[A,B']_{q^{-1}}=0$ (from this follows that $\alpha'$=
$\delta'=0$). See Lemma 2 in reference \cite{aqg1}. Besides 
we require
\begin{equation}
BB'=
\left(
\begin{array}{cc}
\beta\gamma' & 0\\
0 & \gamma\beta'
\end{array}
\right)=q
\left(
\begin{array}{cc}
\beta'\gamma & 0\\
0 & \gamma'\beta
\end{array}
\right)=qBB'.
\label{suem}
\end{equation}
Therefore the following formulas must be fulfilled.
\begin{eqnarray}
b\gamma'=q^{-1}\gamma' a &\;\;\;& a\beta'=q^{-1}\beta' b \label{suemi1}\\
a\beta =q\beta b &\;\;\;& b\gamma=q\gamma a \label{suemi2}\\
\beta\gamma'=q\beta'\gamma  &\;\;\;&\gamma\beta'=q\gamma'\beta. 
\label{suemi3}
\end{eqnarray}

We have two cases. 
\newline
I) For $\epsilon$=$q$
\begin{equation}
\beta=
\left(
\begin{array}{cc}
\beta_{11} & \beta_{12}\\
0 & q\beta_{11}
\end{array}
\right);
\gamma=\beta'=0;
\gamma'=
\left(
\begin{array}{cc}
\gamma'_{11} & \gamma'_{12}\\
0 & q^{-1}\gamma'_{11}
\end{array}
\right).
\end{equation}
Formulas (\ref{suemi1}) and (\ref{suemi2}) follow straightforward, to 
fulfill (\ref{suemi3}) we require
\begin{equation}
a)\;\beta=0,\;\;
\gamma'=
\left(
\begin{array}{cc}
\gamma'_{11} & \gamma'_{12}\\
0 & q^{-1}\gamma'_{11}
\end{array}
\right)\;\mbox{ or }
\end{equation}
\begin{equation}
b)\;\beta=
\left(
\begin{array}{cc}
\beta_{11} & \beta_{12}\\
0 & q\beta_{11}
\end{array}
\right),\;\;
\gamma'=0\;.
\end{equation}
In this case we conclude that either
\begin{equation}  
a)\;\beta=
\left(
\begin{array}{cc}
\beta_{11} & \beta_{12}\\
0 & q\beta_{11}
\end{array}
\right),\;\;\gamma=0,\;\;\beta'=0,\;{\rm and} \;\gamma'=0\;,
\end{equation}
this means $B'=0$ and $BB'=0$, or 
\begin{equation}
b)\;\beta=0,\;\;
\gamma'=
\left(
\begin{array}{cc}
\gamma'_{11} & \gamma'_{12}\\
0 & q^{-1}\gamma'_{11}
\end{array}
\right)
,\;\;\gamma=0\;{\rm and}\;\beta'=0.
\end{equation}
This also means $BB'=0$.\
II) For $\epsilon$=$q^{-1}$
\begin{equation}
\beta'=
\left(
\begin{array}{cc}
\beta'_{11} & \beta'_{12}\\
0 & q^{-1}\beta'_{11}
\end{array}
\right);
\gamma'=\beta=0;
\gamma=
\left(
\begin{array}{cc}
\gamma_{11} & \gamma_{12}\\
0 & q\gamma_{11}
\end{array}
\right)\;.
\end{equation}
Formulas (\ref{suemi1}) and (\ref{suemi2}) follow straightforward, to 
fulfill (\ref{suemi3}) we require
\begin{equation}
a)\;\beta'=0,\;\;
\gamma=
\left(
\begin{array}{cc}
\gamma_{11} & \gamma_{12}\\
0 & q\gamma_{11}
\end{array}
\right)
\end{equation}
\begin{equation}
b)\;\beta'=
\left(
\begin{array}{cc}
\beta'_{11} & \beta'_{12}\\
0 & q^{-1}\beta'_{11}
\end{array}
\right),\;
\gamma=0\;.
\end{equation}
In this case we conclude that either
\begin{equation}
a)\;\beta'=
\left(
\begin{array}{cc}
\beta'_{11} & \beta'_{12}\\
0 & q^{-1}\beta'_{11}
\end{array}
\right),\;\;\gamma=0,\;\;\beta=0,\;{\rm and} \;\gamma'=0\;,
\end{equation}
this means $B=0$ and $BB'=0$, or
\begin{equation}
b)\;\beta=0,\;\;
\gamma=
\left(
\begin{array}{cc}
\gamma_{11} & \gamma_{12}\\
0 & q\gamma_{11}
\end{array}
\right)
,\;\;\gamma'=0\;{\rm and}\;\beta'=0.
\end{equation}
This means $B'=0$ and $BB'=0$.\

Suppose now that one of the matrix $a,b$ is a simplest 
Jordan matrix while the other is a diagonal matrix. A 
conjugation by 
$T=
\left(
\begin{array}{cc}
0  & E\\
E  & 0
\end{array}
\right),$
where $E$ is the identity $2\times 2$ matrix, changes 
$A$=diag$(a,b)$ to diag$(b,a)$, so 
we can suppose that
\begin{equation}
a=\epsilon E+e_{12},\;\;\;b=diag(\mu,1)
\end{equation}
(recall that one of the eigenvalues of $A$ is equal to 
1 and $A$ is an invertible matrix; i.e. $\epsilon, 
\mu\neq 0$).\ 

Firstable, let $\mu\neq q, q^{-1}$, $\delta'=0$, then
\begin{equation}
B=
\left(
\begin{array}{cc}
0  & \beta\\
\gamma  & 0
\end{array}
\right);\;\;
B'=
\left(
\begin{array}{cc}
0  & \beta'\\
\gamma'  & 0
\end{array}
\right)
\end{equation}
We set $BB'=qB'B$ and require formulas (\ref{suemi1})-
(\ref{suemi2}) to hold. From this we obtain the following four 
cases.\newline
I) For $\beta$;
\begin{eqnarray}\nonumber
&A)&  \beta=0\;\;\;\mbox{ provided } \epsilon\neq q\mu, 
\;\;\;\epsilon\neq q,\\ \nonumber
&B) &\beta=
\left(
\begin{array}{cc}
0  & \beta_{12}\\
0  & 0
\end{array}
\right)
\mbox{ provided }
 \epsilon\neq q\mu, \;\;\;\epsilon=q,\\ \nonumber
&C) & \beta=
\left(
\begin{array}{cc}
\beta_{11}  & 0\\
0  & 0
\end{array}
\right)
\mbox{ provided }
 \epsilon=q\mu, \;\;\;\epsilon\neq q,\\ \nonumber
&D) & \beta=
\left(
\begin{array}{cc}
\beta_{11}  & \beta_{12}\\
0  & 0
\end{array}
\right)
\mbox{ provided }
 \epsilon=q\mu, \;\;\;\epsilon=q. \nonumber
\end{eqnarray}
II) For $\gamma$,
\begin{eqnarray}\nonumber
& A) & \gamma=0\;\;\;\mbox{ provided } 
\epsilon\neq q^{-1}\mu, 
\;\;\;\epsilon\neq q^{-1},\\ \nonumber
& B) & \gamma=
\left(
\begin{array}{cc}
0  &  0  \\
0  &  \gamma_{22}
\end{array}
\right)
\mbox{ provided }
 \epsilon\neq q^{-1}\mu, \;\;\;\epsilon=q^{-1},\\ \nonumber
& C) & \gamma=
\left(
\begin{array}{cc}
0  & \gamma_{12}\\
0  & 0
\end{array}
\right)
\mbox{ provided }
 \epsilon=q^{-1}\mu, \;\;\;\epsilon\neq q^{-1},\\ \nonumber
& D) & \gamma=
\left(
\begin{array}{cc}
0  & \gamma_{12}\\
0  & \gamma_{22}
\end{array}
\right)
\mbox{ provided }
 \epsilon=q^{-1}\mu, \;\;\;\epsilon=q^{-1}. \nonumber
\end{eqnarray}
III) For $\beta'$,
\begin{eqnarray}\nonumber
& A) & \beta'=0\;\;\;\mbox{ provided } 
\epsilon\neq q^{-1}\mu, 
\;\;\;\epsilon\neq q^{-1},\\ \nonumber
& B) & \beta'=
\left(
\begin{array}{cc}
0  & \beta'_{12}\\
0  & 0
\end{array}
\right)
\mbox{ provided }
 \epsilon\neq q^{-1}\mu, \;\;\;\epsilon=q,\\ \nonumber
& C) & \beta=
\left(
\begin{array}{cc}
\beta'_{11}  & 0\\
0  & 0
\end{array}
\right)
\mbox{ provided }
\epsilon=q^{-1}\mu, \;\;\;\epsilon\neq q^{-1},\\ \nonumber
& D) & \beta'=
\left(
\begin{array}{cc}
\beta'_{11}  & \beta'_{12}\\
0  & 0
\end{array}
\right)
\mbox{ provided }
 \epsilon=q^{-1}\mu, \;\;\;\epsilon=q^{-1}. \nonumber
\end{eqnarray}
IV) For $\gamma'$,
\begin{eqnarray} \nonumber
& A) & \gamma'=0\;\;\;\mbox{ provided } 
\epsilon\neq q\mu, 
\;\;\;\epsilon\neq q,\\ \nonumber
& B) & \gamma=\left(
\begin{array}{cc}
0  &  0  \\
0  &  \gamma'_{22}
\end{array}
\right)
\mbox{ provided }
 \epsilon\neq q\mu, \;\;\;\epsilon=q,\\ \nonumber
& C) & \gamma=\left(
\begin{array}{cc}
0  & \gamma'_{12}\\
0  & 0
\end{array}
\right)
\mbox{ provided }
 \epsilon=q\mu, \;\;\;\epsilon\neq q.\\ \nonumber
& D) & \gamma\left(
\begin{array}{cc}
0  & \gamma'_{12}\\
0  & \gamma'_{22}
\end{array}
\right)
\mbox{ provided }
 \epsilon=q\mu, \;\;\;\epsilon=q
\end{eqnarray}\nonumber
We can reorganize cases I-IV in the following way.\

i) Let $\epsilon\neq q\mu$, $\epsilon\neq q$. Then we have
\begin{eqnarray}\nonumber
& i.1) &  \beta'_{12}=0 \;\mbox{ or }\; \gamma_{22}=0 \;
\mbox{ for }\; \epsilon\neq q^{-1}\mu; \epsilon=q^{-1};
\epsilon\neq q\mu,\\ \nonumber
& i.2) &  \beta'_{11}=0 \;\mbox{ or }\; 
\gamma_{12}=0 \;\mbox{ for }\; \epsilon= q^{-1}\mu; 
\epsilon\neq q^{-1};\epsilon\neq q\mu,\\ \nonumber
& i.3) &  \gamma_{12}\beta'_{11}=-\beta'_{12}\gamma_{22}
\;\mbox{ for }\;\epsilon=q^{-1}\mu; \epsilon=q^{-1}, \\ \nonumber
& i.4) & \mbox{No extra condition for }\;
\;\epsilon\neq q^{-1}\mu;\; \epsilon\neq q^{-1}
;\;\epsilon\neq q\mu; \epsilon\neq q. \nonumber
\end{eqnarray} 
Studying $i.1)$ we obtain two possible cases.
$$
a)\;\beta'=0\;,\;B'=0\;
$$
$$
b)\;\gamma=0\;,\;B'=0. 
$$
In both cases $BB'$=$0$.\

For $i.2)$ we again obtain two possible cases
$$
a)\;\beta'_{11}=0\;\;\mbox{ then }B'=0. 
$$
$$
b)\;\gamma_{12}=0\;\;\mbox{ then }B=0. 
$$
In both cases $BB'$=$0$.\

For $i.3)$ we obtain,
$$
B=\left(
\begin{array}{cccc}
0  & 0 & 0 & 0\\
0  & 0 & 0 & 0\\
0  & \gamma_{12} & 0 & 0\\
0  & \gamma_{22} & 0 & 0\\
\end{array}
\right)\;\mbox{ and }\;
B'=\left(
\begin{array}{cccc}
0  & 0 & \beta'_{11} & \beta'_{12}\\
0  & 0 & 0 & 0\\
0  & 0 & 0 & 0\\
0  & 0 & 0 & 0\\
\end{array}
\right) 
$$
From this, we conclude again that $BB'$=$B'B$=0.\

In the case $i.4)$ we have $\beta$=$\gamma$=$\beta'$=
$\gamma'$=0, thus $BB'$=0.\

ii) Let $\epsilon\neq q\mu$, $\epsilon=q$.\
In this case we have
$$
B=\left(
\begin{array}{cccc}
0  & 0 & 0 & \beta_{12}\\
0  & 0 & 0 & 0\\
0  & \gamma_{12} & 0 & 0\\
0  & 0 & 0 & 0\\
\end{array}
\right)\;\mbox{ and }\;  
B'=\left(
\begin{array}{cccc}
0  & 0 & \beta'_{11} & 0\\
0  & 0 & 0 & 0\\
0  & 0 & 0 & 0\\
0  & \gamma'_{22} & 0 & 0\\
\end{array}
\right)
$$
From this follows that 
$$
A=diag\left(
\left(
\begin{array}{cc}
q & 1\\
0 & q
\end{array}
\right),
\left(
\begin{array}{cc}
q^2 & 0 \\
0 & 1
\end{array}
\right)\right)\hspace{.5cm}B_1=e_{14},\;B_2=e_{32},\;
B'_1=e_{13}, B'_2=e_{42}
$$
which corresponds to representation (\ref{case3}) in Theorem 1.\

iii) Let $\epsilon\neq q\mu$, $\epsilon\neq q$. In this case we have
$$
B=\left(
\begin{array}{cccc}
0  & 0 & \beta_{11} & 0\\
0  & 0 & 0 & 0\\
0  & 0 & 0 & 0\\
0  & \gamma_{22} & 0 & 0\\
\end{array}
\right)\;\mbox{ and }\;  
B'=\left(
\begin{array}{cccc}
0  & 0 & 0 & \beta'_{12}\\
0  & 0 & 0 & 0\\
0  & \gamma'_{12} & 0 & 0\\
0  & 0 & 0 & 0\\
\end{array}
\right).
$$
From this follows that 
$$
A=diag\left(
\left(
\begin{array}{cc}
q^{-1} & 1\\
0 & q^{-1}
\end{array}
\right),
\left(
\begin{array}{cc}
q^{-2} & 0 \\
0 & 1
\end{array}
\right)\right),\;B_1=e_{13},\;B_2=e_{42},\;
B'_1=e_{14}, B'_2=e_{32}
$$
By applying the maps $q\rightarrow q^{-1}$ and 
$B\rightarrow B'$ we obtain the representation (\ref{case3}) in 
Theorem 1.\

iv). Let $\mu=1$ (namely $\epsilon$=$q\mu$, $\epsilon$=$q$). 
In this case we have\
$$
B=\left(
\begin{array}{cccc}
0  & 0 & \beta_{11} & \beta_{12}\\
0  & 0 & 0 & 0\\
0  & 0 & 0 & 0\\
0  & 0 & 0 & 0\\
\end{array}
\right)\;\mbox{ and }\;  
B'=\left(
\begin{array}{cccc}
0  & 0 & 0 & 0\\
0  & 0 & 0 & 0\\
0  & \gamma'_{12} & 0 & 0\\
0  & \gamma'_{22} & 0 & 0\\
\end{array} 
\right)\;.
$$
From (\ref{suemi3}), follows that $BB'$=0, since $\beta_{11}\gamma'_{12}
+\beta_{12}\gamma'_{22}=0$.\\

Let us now consider the case $\mu=q^{-1}$, we can multiply matrices 
$A$ and $B$ by $q$ and conjugate them by the matrix
$$
diag(1,1,
\left(
\begin{array}{cc}
0 & 1 \\
1 & 0
\end{array}
\right)).
$$
We will obtain an equivalent representation with $\mu=q$.
Thus, it is enough to consider the case
$$
a=\epsilon E+e_{12},\;\;\; b=diag(q,1)
$$
where
$$
\alpha=0,\;\;\; \delta=ce_{12}, \;\;\;c\in {\bf C}.
$$
For 
$\beta=
\left(
\begin{array}{cc}
\beta_{11} & \beta_{12}\\
\beta_{21} & \beta_{22}
\end{array}
\right),$ then we have $\alpha \beta =q\beta b\;$e; 
i.e.
$$
\left(
\begin{array}{cc}
\epsilon & 1 \\
0        & \epsilon
\end{array}
\right)
\left(
\begin{array}{cc}
\beta_{11} & \beta_{12}\\
\beta_{21} & \beta_{22}
\end{array}
\right)=
\left(
\begin{array}{cc}
\epsilon\beta_{11}+\beta_{21}   & \epsilon\beta_{12}
\beta_{22}\\
\epsilon\beta_{21}              & \epsilon\beta_{22}
\end{array}
\right)=
$$
$$
q
\left(
\begin{array}{cc}
\beta_{11}  &  \beta_{12}\\
\beta_{21}  &  \beta_{22}
\end{array}
\right)
\left(
\begin{array}{cc}
q  &  0\\
0  &  1
\end{array}
\right)=q
\left(
\begin{array}{cc}
q\beta_{11}  & \beta_{12}\\
q\beta_{21}  & \beta_{22}
\end{array}
\right),
$$
which implies
\begin{equation}
(\epsilon-q^2)\beta_{11}=-\beta_{21}, 
(\epsilon-q)\beta_{12}=-\beta_{22}\label{sue1}
\end{equation}
\begin{equation}
(\epsilon -q^2)\beta_{21}=0,(\epsilon-q)\beta_{22}=0. \label{sue2}
\end{equation}
If $\epsilon$=$q^2$ then the first equality of (\ref{sue1}) 
gives $\beta_{21}=0$, and if $\epsilon\neq q^2$ then 
the first equality of (\ref{sue2}) gives $\beta_{21}=0$. 
Therefore $\beta_{21}=0$ in any case. In the same way 
$\beta_{22}=0$ and (\ref{sue1}), (\ref{sue2}) are equivalent to
\begin{equation}
(\epsilon-q^2)\beta_{11}=0,\;(\epsilon-q)\beta_{12}=0, \label{su1}
\end{equation}
\begin{equation}
\beta_{21}=0,\; \beta_{22}=0. \label{su2}
\end{equation}
Analogously for the matrix 
$\gamma=
\left(
\begin{array}{cc}
\gamma_{11} & \gamma_{12} \\
\gamma_{21} & \gamma_{22}
\end{array}
\right)$
we have $b\gamma$=$q\gamma a$; i.e.
$$
\left(
\begin{array}{cc}
q & 0\\
0 & 1
\end{array}
\right)
\left(
\begin{array}{cc}
\gamma_{11} & \gamma_{12} \\
\gamma_{21} & \gamma_{22}
\end{array}
\right)=
\left(
\begin{array}{cc}
q\gamma_{11} & q\gamma_{12} \\
\gamma_{21} & \gamma_{22}
\end{array}
\right)=
$$
$$
q
\left(
\begin{array}{cc}
\gamma_{11} & \gamma_{12} \\
\gamma_{21} & \gamma_{22}
\end{array}
\right)
\left(
\begin{array}{cc}
\epsilon & 1 \\
0 & \epsilon
\end{array}
\right)=q
\left(
\begin{array}{cc}
\epsilon\gamma_{11} & \gamma_{11}+\epsilon\gamma_{12} \\
\epsilon\gamma_{21} & \gamma_{21}+\epsilon\gamma_{22}
\end{array}
\right).
$$
This implies
\begin{equation}
q(1-\epsilon)\gamma_{11}=0,\;q(1-\epsilon)
\gamma_{12}=q\gamma_{11}\label{suem1}
\end{equation}
\begin{equation}
(1-q\epsilon)\gamma_{21}=0,\;(1-q\epsilon)
\gamma_{22}=q\gamma_{21}.\label{suem2}
\end{equation}
Again, if $\epsilon=1$ then by the second equality of 
(\ref{suem1}), $\gamma_{11}=0$ and if $\epsilon\neq 1$ then by 
the first one $\gamma_{11}=0$. In the same way 
$\gamma_{21}=0$ and (\ref{suem1}),(\ref{suem2}) are equivalent to
\begin{equation}
\gamma_{11}=0,\;\;(1-\epsilon)\gamma_{12}=0, \label{kha1}
\end{equation}
\begin{equation}
\gamma_{21}=0,\;\;(1-q\epsilon)\gamma_{22}=0.\label{kha2}
\end{equation}
Now if $\epsilon\neq q^{-1},1,q,q^2$ then by (\ref{su1}),
(\ref{su2}) and (\ref{kha1}), (\ref{kha2}) $\beta=\gamma=0$ and the 
representation has the form
\begin{equation}
A=diag\left(
\left(
\begin{array}{cc}
\epsilon & 1 \\
0	 & \epsilon\\
\end{array}
\right), q, 1\right),\;\;\; B=e_{34}.
\end{equation}
In this case $\left(B(A)\cdot B'(A)\right)\cap $
$\left(B(A)'\cdot B(A)\right)=0$; namely the representation is not
admissible.\

Finally, let us consider four last possibilities.\

1. $\epsilon$=$q^{-1}$. By (\ref{su1}) and  (\ref{su2}) we have $\beta=0$
 and by (\ref{kha1}) and (\ref{kha2}), $\gamma=ce_{22}$. From this follows
that
$$
A=diag\left(
\left(
\begin{array}{cc}
q^{-1} & 1 \\
0      & q^{-1}\\
\end{array}
\right),q,1\right), \;\;\;B_1=e_{42} , \;\;\;B_2=e_{34}.
$$
If we multiply $A$ by $q$ and conjugate it by $T$=
diag$(1, q^{-1},1,1)$ we will obtain an equivalent 
representation $A$=$diag(1,1,q^2,q)+e_{12}$, $B_1$=
$e_{42}$, $B_{2}$=$e_{34}$. Using conjugations by 
matrices $E-e_{ii}-e_{jj}+e_{ij}+e_{ji}$ we can change 
indices with the help of permutation $1\rightarrow 3$, 
$2\rightarrow 4$, $3\rightarrow 1$, $4\rightarrow 2$. 
Therefore, $e_{42}\rightarrow e_{24}$, 
$e_{34}\rightarrow e_{12}$ and we have the 
representation $A$=$diag(q^{2},q,1,1)$+$e_{34}$,
$B_1=e_{24}$, $B_2=e_{12}$, $B'_1=e_{21}$, 
$B'_2=e_{32}$. In this case $\left(B(A)\cdot B'(A)\right)\cap $
$\left(B(A)'\cdot B(A)\right)=0$; namely the representation is not
admissible.\

2. $\epsilon=1$. By (\ref{su1}) and (\ref{su2}) we again have 
$\beta=0$ and by (\ref{kha1}) and (\ref{kha2}), $\gamma=Ce_{12}$. From 
this the representation has the following form.
$$
A=diag\left(
\left(
\begin{array}{cc}
1 & 1 \\
0 & 1\\
\end{array}
\right),q,1\right), \;\;\;B=e_{32} , 
\;\;\;B_2=e_{34}.
$$ 
and therefore $B(A)^2=0$.\

3. $\epsilon=q$. By (\ref{kha1}) and (\ref{kha2}) we have 
$\gamma=0$ and by (\ref{su1}) and (\ref{su2}) $\beta=ce_{12}$. From this 
the representation has the form
$$
A=diag\left(
\left(
\begin{array}{cc}
q & 1 \\
0 & q\\
\end{array}
\right),q,1\right), \;\;\;B_1=e_{14} , 
\;\;\;B_2=e_{34}
$$
and again $B(A)^2=0$.\

4. $\epsilon=q^2$. By (\ref{kha1}) and (\ref{kha2}) we have 
$\gamma=0$ and equalities (\ref{su1}) and (\ref{su2}) imply 
$\beta=ce_{11}$. So the representation has the form
$$
A=diag\left(
\left(
\begin{array}{cc}
q^2 & 1 \\
0 & q^2\\
\end{array}
\right),q,1\right), \;\;\;B_1=e_{13} , 
\;\;\;B_2=e_{34}.
$$
In this case 
$\left(B(A)\cdot B'(A)\right)\cap $
$\left(B(A)'\cdot B(A)\right)=0$; namely the representation
is not admissible.
\hfill$\Box$.\
\section {$GL_2$ representations.}
\begin{theorem}$\!\!\!.$
Each irreducible finite dimensional algebra 
representation of the quantum $GL_2$, $q^m\neq 1,$ is one 
dimensional.
\end{theorem}
{\it Proof}. Let $c_{ij}\rightarrow C_{ij}$ be a finite 
dimensional irreducible representation of the quantum 
$GL_2$, where $C_{ij}$ are $n\times n$ matrices acting 
on the $n$-dimensional space $V$. This means that the 
matrices $C_{ij}$ satisfy the relations of $GL_2$: 
\begin{equation}
C_{11}C_{12}=C_{12}C_{11},\;\; C_{21}C_{11}=qC_{11}C_{21},\label{theo21}
\end{equation}
\begin{equation}
C_{22}C_{12}=qC_{12}C_{22},\;\; C_{21}C_{22}=C_{22}C_{21},\label{theo22}
\end{equation}
\begin{equation}
C_{21}C_{12}=qC_{12}C_{21},\;\;C_{22}C_{11}-C_{11}C_{22}
=(q-1)C_{12}C_{21}\label{theo23}
\end{equation}
and the matrix $det_q=C_{11}C_{22}-C_{12}C_{21}$ is 
invertible.\

From the relations (\ref{theo21})-(\ref{theo23}) follow that $C_{12}V$ 
is an invariant subspace:
\begin{equation}
C_{11}(C_{12}V)=C_{12}(C_{11}V)\subseteq C_{12}V,
\end{equation}
\begin{equation}
C_{22}(C_{12}V)=qC_{12}C_{22}V=C_{12}(qC_{22}V)
\subseteq C_{12}V,
\end{equation}
\begin{equation}
C_{21}(C_{12}V)=qC_{12}(C_{21}V)=C_{12}(qC_{21}V)
\subseteq C_{12}V.
\end{equation}
Therefore either $C_{12}=0$ or $C_{12}$ is an 
invertible matrix. In the same way either $C_{21}=0$ or 
$C_{21}$ is invertible.\

If both matrices $C_{12}$, $C_{21}$ are equal to zero, then 
the matrices $C_{11}$, $C_{22}$ commute therefore they have 
a common eigenvector $v$ and $Cv$ is an invariant subspace, 
so $Cv=V$, $dim V=1$.\

Suppose that $C_{21}=0$ and $C_{12}$ is invertible. Then 
$det_q=C_{11}C_{22}$ and both matrices $C_{11}$, $C_{22}$ 
are invertible. Now $x\rightarrow C_{22}$, 
$y\rightarrow C_{12}$ is a representation of the $q$-spinor 
with invertible matrices which is a contradiction. Recall that
if $q^m\neq 1$, and one of the matrices in the $q$-spinor is
invertible; then the second one must be nilpotent.\

Suppose that $C_{12}=0$ and $C_{21}$ is invertible. Then 
$det_q=C_{11}C_{22}$ and both matrices $C_{11}$, $C_{22}$ 
are invertible. Now $x\rightarrow C_{11}$, 
$y\rightarrow C_{21}$ is a representation of the $q^{-1}$-
spinor with invertible matrices which is a contradiction.\

Finally, let $C_{12}$, $C_{21}$ be invertible matrices and 
$C_{11}$, $C_{22}$ be nilpotent ones. We have the following 
relation
\begin{equation}
[C_{11},C_{22}]=(q-1)C_{12}C_{21}=\epsilon,
\end{equation}
here $\epsilon$ is an invertible matrix, such that
\begin{equation}
\epsilon C_{11}=qC_{11}\epsilon .
\end{equation}
Using this relation and induction by $k$ we can prove that,
\begin{equation}
[C^k_{11}, C_{22}]=q^{[k]}C^{k-1}_{11}\epsilon, \label{kmutator}
\end{equation}
where
\begin{equation}
x^{[k]}=1+x+...+x^{k-1}=\frac{x^{k}-1}{x-1}=
x^{[k-1]}\cdot x +1.
\end{equation}

Indeed, if $k$ is the smallest number such that $C^k_{11}=0$, 
then (\ref{kmutator}) gives a contradiction: $[C^k_{11}, C_{22}]$=
$q^{[k]}C^{k-1}_{11}\epsilon=0$ and $C^{k-1}_{11}=0$ since 
$q^{[k]}$=$\frac{1-q^{k}}{1-q^{-2}}\neq 0$ and $\epsilon$ is 
invertible.\hfill$\Box$.\

From this follows straightforward that every finite 
dimensional representation of the quantum 
$GL_2, q^m\neq 1,$ is triangular; i.e.  it is equivalent to
 a representation by triangular matrices $c_{ij}
\rightarrow C_{ij}$.\

\begin{corollary}$\!\!\!.$
 For every finite dimensional representation 
$c_{ij}\rightarrow C_{ij}$ of the quantum 
$GL_2, q^m\neq 1$, the elements $C_{11}$, $C_{22}$ are 
invertible, while $C_{12}$, $C_{21}$ are nilpotent.
\end{corollary}
{\it Proof}. We can suppose that $C_{ij}$ 
are triangular matrices. In this case 
the matrix
\begin{equation}
(1-q)^{-1}(C_{11}C_{22}-C_{22}C_{11})
\end{equation}
has only zero entries on the main diagonal. This matrix is 
equal to $C_{12}C_{21}$. From this follows that the main 
diagonal of $C_{11}C_{22}$ and that of the invertible 
matrix $det_q$=$C_{11}C_{22}-C_{12}C_{21}$ coincide. This 
means that $C_{11}$ and $C_{22}$ have no zero terms on the 
main diagonal and therefore they are invertible.
\hfill$\Box$\
\begin{theorem}$\!\!\!.$ Let $c_{ij}
\rightarrow C_{ij}$ and $c_{ij}\rightarrow C'_{ij}$ be 
two representations of $GL_2$ in ${\it C}(1,3)$. Then Hopf 
algebra actions
\begin{equation}
c_{ij}\cdot v=\sum_kC_{ik}vC^*_{kj}
\end{equation}
and
\begin{equation}
c_{ij}*v=\sum_kC'_{ik}vC'^*_{kj}
\end{equation}
are equivalent if and only iff
$$C'_{11}=uC_{11}u^{-1}\alpha_1,\;
C'_{12}=uC_{12}u^{-1}\alpha_2,\;$$
\begin{equation}
C'_{21}=uC_{21}u^{-1}\alpha_1,\;
C'_{22}=uC_{22}u^{-1}\alpha_2,
\end{equation}
for some nonzero complex numbers $\alpha_1$, $\alpha_2$
 and invertible $u\in {\it C}(1,3)$. If $u=1,$ then 
the actions coincide.
\end{theorem}
{\it Proof}. The proof follows like in Theorem 2, reference \cite{aqg1}.\

In terms of modules, this result says that the equivalence 
of representations means that corresponding modules 
$V_1$, $V_2$ are related by formula $V_1\simeq V_2
\otimes {\it U}$, where ${\it U}$ is any one dimensional 
module.\

\section{Invariants and the operator algebra.}

For a given representation $c_{ij}\rightarrow C_{ij}$ we denote by
$\Re $ an {\it operator algebra} i.e. a 
subalgebra of $C(1,3)$ generated by $C_{ij}.$ 
Recall that the algebra of invariants of an action is defined in the
following way
\begin{equation}
Inv=\{ v\in C | \forall h\in H \ \ \ h\cdot v=\varepsilon (h)v\} .
\label{inv}
\end{equation}
being $H$ any Hopf algebra and $\epsilon (h)$ the corresponding
counit. On the other hand the Invariant algebra equals the
 centralizer of $\Re$ in ${\it C}(1,3)$.\

In this Section we present five ingredients for 
every representation of the quantum $GL_2$ by Dipper-Donkin with 
{\it nonzero perturbation}: 
the values of $C_{ij},$ the 
matrix form of the operator algebra $\Re $, its 
dimension, the invariants of the inner action 
defined by this representation $I$, and the value of
the quantum determinant.\

To obtain the full classification presented in reference \cite{math}, 
from where we extract the representations given in this Section, 
Theorem 1 and Figure 2 are used. Additional information
 (i.e. $\Re$, $I$, etc) are derived from 
Theorem 2 and Corollary 1. Theorem 3 is intended to address the 
question about minimal nonequivalent representation for $GL_2$ by
Dipper and Donkin on ${\it C}(1,3)$, this question remains open, so far.\

{\bf CASE 1).}\\

$\begin{array}{l}
         d=diag(q^2,q,q,1)\\ C_{12}=qe_{13}-\mu e_{24}\\
C_{21}=-\mu e_{21}+e_{43}\\ 
C_{11}=diag(1,q^{-1}, 1, q^{-1})\\
C_{22}=diag(q^2,q^2,q,q)-q\mu e_{23}\;\;\;\;\end{array}$
$;\;\;\;\Re \cong\left(\matrix{*&*\cr 0&*\cr}\right)\otimes
\left(\matrix{*&*\cr 0&*\cr}\right);$
\vspace*{1cm}
\newline
\hspace*{4cm}$dim \Re  =9\;\;$$;\;\;\;\;
$${\rm Invariants}\cong C.$\\

This corresponds to CASE 4) in Theorem 1.\\ 

{\bf CASE 2).}\\

$\begin{array}{l}
         d=diag(q^2,q,q,1)\\ 
C_{12}=qe_{12}+\mu e_{34}\\
C_{21}=\mu e_{31}+e_{42}\\
C_{11}=diag(1,1, q^{-1}, q^{-1})\\
C_{22}=diag(q^2,q,q^2,q)+q\mu e_{32}\;\;\;\;\end{array}$
$;\;\;\;\Re \cong\left(\matrix{*&*\cr 0&*\cr}\right)\otimes
\left(\matrix{*&*\cr 0&*\cr}\right);$
\vspace*{1cm}
\newline
\hspace*{4cm}$dim \Re  =9\;\;$$;\;\;\;\;
$${\rm Invariants}\cong C.$\\

This corresponds to CASE 5) in Theorem 1.\\

For CASE 6), in Theorem 1, we find that there exist no set
$\{c_{ij},d\}$, $1\leq i,j \leq 2$, that fulfills the algebra
in Figure 1.\

All the possible representations for the quantum $GL_2$ by Dipper-
Donkin are reported elsewhere \cite{math}. From there, we can deduce
the following.\

a) Only for nonzero perturbation representations of $GL_2$, 
$dim \Re =9$; which turns out to be the maximal possible dimension 
of $\Re$ for the action of $GL_2$ on ${\it C}(1,3)$.\

b) The maximal dimension for $I$ is 6.\

c) There is only one, zero perturbation, possible case in the set 
of all representations of $GL_2$ by Dipper-Donkin, on ${\it C}(1,3)$, 
for which $dim \Re$=8. This is as follows
$$\begin{array}{ll}
d=diag(q^2,q,1,1) & C_{12}=\alpha e_{12}+\beta e_{23}+\gamma e_{24}\\
C_{21}=0 & C_{11}={\bf 1} \\
C_{22}=q^2e_{11}+qe_{22}+e_{33}+e_{44}. & 
\end{array}$$\

Besides, for this case we know that\

$$\Re =\left(
\begin{array}{cccc}
* & * & * & * \\
0 & * & * & * \\
0 & 0 & \epsilon & 0 \\
0 & 0 & 0 & \epsilon 
\end{array}
\right)\;;\;\;I=C.$$\

d) There is only one, zero perturbation, possible case in the set of
all representations of $GL_2$ by Dipper-Donkin, on ${\it C}(1,3)$, for
which $dim\Re $=3. This is as follows
$$\begin{array}{ll}
d=diag(q^2,q,q,1) & C_{12}=0\\
C_{21}=0 & C_{11}=e_{11}+\alpha_{1}e_{22}+m\alpha_{2}e_{33}+
\alpha_3 e_{44} \\
C_{22}=q^2e_{11}+q\alpha^{-1}_{1}e_{22}+\alpha^{-1}_{2}e_{33}+
{\alpha}^{-1}_{3}e_{44}. & 
\end{array}$$\

Besides, for this case we know that\

$$\Re =\left(
\begin{array}{cccc}
* & 0 & 0 & 0\\
0 & \epsilon & 0 & 0 \\
0 & 0 & \epsilon & 0 \\
0 & 0 & 0 & *
\end{array}
\right)\;;\;\;
I=
\left(
\begin{array}{cccc}
\alpha & 0 & 0 & 0\\
0 & \beta & \gamma & 0 \\
0 & \delta & \epsilon & 0 \\
0 & 0 & 0 & *
\end{array}
\right)
.$$\ 

%Incidentaly, this is the minimum possible dimension for $\Re$.\

e) For the representations wherein $d$=$diag(\alpha,q^2,q,1)$,
$\alpha\neq 0$, $q^{-1}$, 1, $q$, $q^2$, $q^3$, which correspond to
CASE 5) in Theorem (1) \cite{aqg1}, always $dim \Re$=$6$ and 
$I=C\oplus C$.\

f) For the representations wherein $d=diag(q^3,q^2,q,1)$ which 
correspond to CASE 4) in Theorem 1 \cite{aqg1}, always $dim \Re$=$7$ and
$I=C$.\

{\rm \section{Acknowledgments.}

The author wishes to thank CONACYT for partial support under grant
4336-E, and Vladislav Kharchenko for helpful discussion.

}

\end{document}